\documentclass[12pt,twoside,reqno]{amsart}
\usepackage[colorlinks=false,citecolor=blue]{hyperref}
\usepackage{mathptmx, amsmath, amssymb, amsfonts, amsthm, mathptmx, enumerate, color,mathrsfs}
\usepackage{graphicx}

\usepackage{epstopdf}
\newtheorem{theorem}{Theorem}[section]

\newtheorem{corollary}{Corollary}[section]

\numberwithin{equation}{section}

\begin{document}
	
\title{Some functional identities characterizing two-sided centralizers and two-sided generalized derivations on triangular algebras}
\author{Amin Hosseini}

\begin{abstract}
Let $\mathcal{T}$ be a unital triangular algebra, let $n > 1$ be an integer, let $\gamma$ be an invertible element of $Z(\mathcal{T})$, the center of $\mathcal{T}$, and let $\Psi, \Omega:\mathcal{T}\rightarrow \mathcal{T}$ be additive mappings satisfying \begin{align*}
\Psi(X^n) = \gamma X^{n - 1}\Omega(X) = \gamma \Omega(X) X^{n - 1}\end{align*}
for all $X \in \mathcal{T}$. If $\Omega(\textbf{1}) \in Z(\mathcal{T})$, then $\Psi$ and $\Omega$ are two-sided centralizers on $\mathcal{T}$ and also $\Psi = \gamma \Omega$. Moreover, using a functional identity, a characterization of two-sided generalized derivations is presented. Some other related results are also discussed.
\end{abstract}

\maketitle

\pagestyle{myheadings}
\markboth{\centerline {}}
{\centerline {}}
\bigskip
\bigskip

\section{Introduction and Preliminaries}

Let $\mathcal{R}$ be a ring. An additive mapping $T:\mathcal{R} \rightarrow \mathcal{R}$ is said to be a \emph{left} (resp. \emph{right}) \emph{centralizer} or \emph{multiplier} of $\mathcal{R}$ if $T(xy) = T(x)y$ (resp. $T(xy) = x T(y)$) for all $x, y \in \mathcal{R}$. We call $T$ a \emph{two-sided centralizer} (some authors call it a centralizer) if it is both left and right centralizer. An additive mapping $T$ on $\mathcal{R}$ is called a \emph{Jordan left} (resp. \emph{right}) centralizer if $T(x^2) = T(x) x$ (resp. $T(x^2) = x T(x)$) for all $x \in \mathcal{R}$, and it is called a \emph{Jordan two-sided centralizer} if $T(x^2) = T(x)x = x T(x)$ for all $x \in \mathcal{R}$. An additive mapping $T$ on $\mathcal{R}$ is called a \emph{Jordan centralizer} if $T(x \circ y) = T(x) \circ y = x \circ T(y)$ for all $x, y \in \mathcal{R}$, where $x \circ y = xy + yx$. Since the product $\circ$ is commutative, there is no difference between the left and right Jordan centralizers. Obviously, every left (resp. right) centralizer is a Jordan left (resp. right) centralizer. But the converse is, in general, not true.

Centralizers are very important both in theory and applications and these mappings have been studied in the general framework of prime rings and semiprime rings. Recall that a ring $\mathcal{R}$ is called prime if for $x, y \in \mathcal{R}$, $x \mathcal{R} y = \{0\}$ implies that $x = 0$ or $y = 0$, and is semiprime if for $x \in \mathcal{R}$, $x \mathcal{R} x = \{0\}$ implies that $x = 0$. A ring $\mathcal{R}$ is said to be $n$-torsion free, where $n > 1$ is an integer, if for $x \in \mathcal{R}$, $nx = 0$ implies that $x = 0$. Similarly, a module $\mathcal{M}$ is said to be $n$-torsion free, where $n > 1$ is an integer, if for $\mathfrak{m} \in \mathcal{M}$, $n\mathfrak{m} = 0$ implies that $\mathfrak{m} = 0$.

As a pioneering work, Zalar \cite{Z} proved that if $\mathcal{R}$ is a 2-torsion free semiprime ring and $T$ is a Jordan centralizer on $\mathcal{R}$, then $T$ is a two-sided centralizer. After that, this question under what conditions a map becomes a two-sided centralizer attracted much attention of mathematicians. So far, many mathematicians have extensively investigated the characterization of centralizers using functional identities on rings as well as algebras and this notion has been studied from other perspectives as well. To read more about centralizers, we refer the readers to some recent papers \cite{A, C, F1, F*, Fo1, Fo2, G, H, K1, K2, L, V}, where further references can be found. Jordan centralizers on triangular rings (or algebras) have also been investigated. Some results concerning Jordan centralizers can be found in \cite{As, G1, H2, H, L1, L2, F2}, and the references therein.

As the main tool in this paper, we apply the theory of functional identities. The theory of functional identities considers set-theoretic mappings on rings (or algebras) that satisfy some identical relations. When treating such relations one usually concludes that the form of the mappings involved can be described, unless the ring is very special. We refer the reader to \cite{br1, br2} for a full treatment of this theory.

In the present paper, we characterize two-sided centralizers and two-sided generalized derivations on triangular rings via some functional identities. Let $\mathcal{A}$ and $\mathcal{B}$ be two algebras and let $\mathcal{M}$ be an $(\mathcal{A}, \mathcal{B})$-bimodule which is faithful as a left $\mathcal{A}$-module as well as a right $\mathcal{B}$-module. The ring
\begin{align*}
\mathcal{T} = Tri(\mathcal{A}, \mathcal{M}, \mathcal{B}) := \left \{\left [\begin{array}{cc}
a & m\\
0 & b
\end{array}\right ] \ : \ a \in \mathcal{A}, \ b \in \mathcal{B}, \ m \in \mathcal{M} \right\}
\end{align*}
under the usual matrix operations is said to be a \emph{triangular algebra}. It is clear that $\mathcal{T} = Tri(\mathcal{A}, \mathcal{M}, \mathcal{B})$ is unital if and only if both $\mathcal{A}$ and $\mathcal{B}$ are unital. Moreover, the centre $Z(\mathcal{T})$ of $\mathfrak{T} = Tri(\mathcal{A}, \mathcal{M}, \mathcal{B})$ is of the form below.
\begin{align*}
Z(\mathcal{T}) = \big\{a \oplus b \ : \ a \in Z(\mathcal{A}), b \in Z(\mathcal{B}), am = mb \ for \ all \ m \in \mathcal{M} \big\}.
\end{align*}
For more details in this regard, see \cite{C2, C3}. Let $\mathcal{A}$ and $\mathcal{B}$ be normed algebras and let $\mathcal{M}$ be a normed $(\mathcal{A}, \mathcal{B})$-bimodule. We can define a norm on the triangular algebra $\mathcal{T} = Tri(\mathcal{A}, \mathcal{M}, \mathcal{B})$ as follows:
\begin{align*}
\Big\|\left [\begin{array}{cc}
a & m\\
0 & b
\end{array}\right ]\Big\| = \|a\|_{\mathcal{A}} + \|m\|_{\mathcal{M}} + \|b\|_{\mathcal{B}}.
\end{align*}
A routine calculation shows that $\mathcal{T}$ is in fact a normed algebra. If $\mathcal{A}$ and $\mathcal{B}$ are Banach algebras and $\mathcal{M}$ is a Banach $(\mathcal{A}, \mathcal{B})$-bimodule, then the triangular algebra $\mathcal{T} = Tri(\mathcal{A}, \mathcal{M}, \mathcal{B})$ is a Banach algebra with respect to the abovementioned norm. We shall call such an algebra a triangular Banach algebra. For more details, see, e.g. \cite{F2}.

Following \cite{F} and \cite[Lemma 1.5]{J}, we say that a ring $\mathcal{R}$ satisfies Condition (P) if $x r x = 0$ for all $x \in \mathcal{R}$ and for some $r \in \mathcal{R}$, then $r = 0$. It is clear that every unital ring satisfies Condition (P). Hosseini \cite{H1} has also proved that every semiprime ring satisfies Condition (P). But, the converse is not true in general. In deed, according to \cite[Example 3]{F}, there exists a non-semiprime ring which satisfies Condition (P).

In the following, we explain the main results of this article.
Let $\mathcal{A}$ and $\mathcal{B}$ be unital algebras and let $\mathcal{M}$ be a faithful $(\mathcal{A}, \mathcal{B})$-bimodule. Let $n > 1$ be an integer, let $\gamma$ be an invertible element of $Z(\mathcal{T})$ and let $\Psi, \Omega:\mathcal{T}\rightarrow \mathcal{T}$ be additive maps satisfying \begin{align*}
\Psi(X^n) = \gamma X^{n - 1}\Omega(X) = \gamma \Omega(X) X^{n - 1}\end{align*}
for all $X \in \mathcal{T}$. If $\Omega(\textbf{1}) \in Z(\mathcal{T})$, then $\Psi$ and $\Omega$ are two-sided centralizers on $\mathcal{T}$ and also $\Psi = \gamma \Omega$. As a consequence of this result, we prove that under the stated conditions, if $\Psi, \Omega:\mathcal{T} \rightarrow \mathcal{T}$ are additive maps satisfying \begin{align*}
\Psi(X^n) = \gamma X\Omega(X^{n-1}) = \gamma \Omega(X^{n - 1}) X\end{align*}
for all $X \in \mathcal{T}$ and further $\Omega(\textbf{1}) \in Z(\mathcal{T})$, then $\Psi$ and $\Omega$ are two-sided centralizers on $\mathcal{T}$ and also $\Psi = \gamma \Omega$. These results lead to conclusions about the automatic continuity of two-sided centralizers on triangular normed algebras.

Another main result of this study is to present a functional identity characterizing two-sided generalized derivations on triangular rings. Let $\mathcal{R}$ be a ring. Recall that an additive map $F : \mathcal{R} \rightarrow \mathcal{R}$  is called an \emph{l-generalized derivation} if there exists a derivation $d : \mathcal{R} \rightarrow \mathcal{R}$ such that $F(xy) = F(x) y + x d(y)$ for all $x, y \in \mathcal{R}$, which is the notion introduced by Bre$\check{s}$ar \cite{br}. Also, an additive map $G : \mathcal{R} \rightarrow \mathcal{R}$ is called an \emph{r-generalized derivation} associated with a derivation $g : \mathcal{R} \rightarrow \mathcal{R}$, if $G(xy) = g(x) y + x G(y)$ for all $x, y \in \mathcal{R}$. We say that a linear map $F$ is a \emph{two-sided generalized derivation} if it is both an l-generalized derivation associated with a derivation $d_1$ and an r-generalized derivation associated with a derivation $d_2$. For more material about two-sided generalized derivations, see, e.g. \cite{Ab, G2, H}. It is routine to see that if $F$ is a left (resp. right) generalized derivation associated with a derivation $d$, then $T = F - d$ is a left (resp. right) centralizer. Our result about the characterization of two-sided generalized derivations reads as follows. Let $n > 1$ be an integer and let $\Psi, \Omega:\mathcal{T}\rightarrow \mathcal{T}$ be linear mappings satisfying \begin{align*}
2\Psi(X^n) =  X^{n - 1}\Omega(X) + \Omega(X) X^{n - 1}\end{align*}
for all $X \in \mathcal{T}$. If $\Omega(\textbf{1}) \in Z(\mathcal{T})$, then $\Psi$ and $\Omega$ are two-sided generalized derivations on $\mathcal{T}$. In particular, if $\Psi = \Omega$, then $\Omega (= \Psi)$ is a two-sided centralizer.

\section{Results and Proofs}

Inspired by \cite[Theorem 2.7]{Hoss}, we present the theorem below.

\begin{theorem}\label{1}
Let $n > 1$ be an integer, let $\gamma$ be an invertible element of $Z(\mathcal{T})$ and let $\Psi, \Omega:\mathcal{T}\rightarrow \mathcal{T}$ be additive mappings satisfying \begin{align*}
\Psi(X^n) = \gamma X^{n - 1}\Omega(X) = \gamma \Omega(X) X^{n - 1}\end{align*}
for all $X \in \mathcal{T}$. If $\Omega(\textbf{1}) \in Z(\mathcal{T})$, then $\Psi$ and $\Omega$ are two-sided centralizers on $\mathcal{T}$ and also $\Psi = \gamma \Omega$.
\end{theorem}
\begin{proof}According to the aforementioned assumption, we have
\begin{align}
\Psi(X^n) = \gamma X^{n - 1}\Omega(X),\ for \ all \ X \in \mathcal{T}.
\end{align}
Let $C$ be an arbitrary element of $Z(\mathcal{T})$. Replacing $X$ by $X + C$ in (2.1), we have
\begin{align*}
\Psi\Big(\sum_{k = 0}^{n}\Big(^{n}_{k}\Big)X^{n - k}C^{k}\Big)  = \gamma \sum_{k = 0}^{n - 1}\Big(^{n - 1}_{k}\Big)X^{n - k - 1}C^{k}\left(\Omega(X) + \Omega(C)\right)
\end{align*}
Using equation (2.1) and collecting together terms of the above-mentioned expressions involving the same number of factors of $C$, it is obtained that
\begin{align}
\sum_{k = 1}^{n - 1}\Phi_k(X, C) = 0, \ for \ all \ X \in \mathcal{T},
\end{align}
where
\begin{align*}
\Phi_k(X, C) = \Big(^{n}_{k}\Big)\Psi(X^{n - k}C^{k}\Big) - \gamma \Big(^{n-1}_{k}\Big)C^{k}X^{n - k - 1}\Omega(X) - \gamma \Big(^{n-1}_{k-1}\Big)C^{k - 1}X^{n - k}\Omega(C).
\end{align*}
Replacing $C$ by $C, 2C, ..., (n - 1)C$ in equation (2.2), we obtain a system of $(n - 1)$ homogeneous equations as follows:
\[
\left\lbrace
  \begin{array}{c l}
     \sum_{i = 1}^{n- 1}\Phi_{i}(X,C) = 0,\\
\\
 \sum_{i = 1}^{n- 1}\Phi_{i}(X,2 C) = 0,\\  & 
  \\
  .\\
  .\\
  .\\
  \sum_{i = 1}^{n - 1}\Phi_{i}(X,(n - 1) C) = 0.
  \end{array}
 \right. \]

We see that the coefficient matrix of the above system is:\\
\begin{align*}
Y = \left [\begin{array}{ccccccc}
1 & 1 & 1 & . & . & . & 1\\
2 & 2^{2} & 2^{3} & . & . & . & 2^{n - 1}\\
. & . & . & & . & & .\\
. & . & . & & . & & .\\
. & . & . & & . & & .\\
(n - 1) & (n - 1)^{2} & (n - 1)^{3} & . & . & . & (n - 1)^{n -1}
\end{array}\right ]
\end{align*}
We know that the determinant of a square Vandermonde matrix is nonzero, and it implies that the above-mentioned system has only a trivial solution. In particular, $\Phi_{n - 1}(X,C) = 0$ and $\Phi_{n - 2}(X,C) = 0$. Since $\Phi_{n - 1}(X,C) = 0$, we have
\begin{align*}
0 = & \Phi_{n - 1}(X,C) = \Big(_{n - 1}^{n}\Big)\Psi(X C^{n - 1}) - \gamma C^{n - 1} \Omega(X) - \gamma \Big(_{n - 2}^{n-1}\Big)C^{n - 2}X \Omega(C)
\end{align*}
for all $X \in \mathcal{T}$. Putting $C = \textbf{1}$ in the above equation, we get that
\begin{align}
\Psi(X) = \frac{ \gamma}{n} \Omega(X) + \frac{\gamma (n - 1)}{n} X \Omega(\textbf{1}), \ for \ all \ X \in \mathcal{T}.
\end{align}
It follows from (2.3) that
\begin{align}
\Psi(X^2) = \frac{ \gamma}{n} \Omega(X^2) + \frac{\gamma (n - 1)}{n} X^2 \Omega(\textbf{1}), \ for \ all \ X \in \mathcal{T}.
\end{align}

On the other hand, since $\Phi_{n - 2}(X,C) = 0$, we have the following expressions:
\begin{align*}
0 = & \Phi_{n - 2}(X,C) = \Big(_{n - 2}^{n}\Big)\Psi(X^{2}C^{n - 2}) - \gamma \Big(_{n - 2}^{n-1}\Big)C^{n - 2}X \Omega(X) - \gamma \Big(_{n - 3}^{n-1}\Big)C^{n - 3}X^{2} \Omega(C)
\end{align*}
for all $X \in \mathcal{T}$. Putting $C = \textbf{1}$ in the above equation, we get that
\begin{align}
\Psi(X^{2}) = \frac{2 \gamma}{n}X \Omega(X) + \frac{\gamma (n - 2)}{n} X^{2} \Omega(\textbf{1}), \ for \ all \ X \in \mathcal{T}.
\end{align}
Comparing (2.4) and (2.5) and using the assumption that $\gamma$ is invertible, we get that
\begin{align}
\Omega(X^{2}) = 2 X \Omega(X) - \Omega(\textbf{1}) X^2 , \ for \ all \ X \in \mathcal{T}.
\end{align}
Considering the functional identity $\Psi(X^n) = \gamma \Omega(X^{n-1})X$ for all $X \in \mathcal{T}$, and using an argument similar to the one above, we deduce that
\begin{align}
\Omega(X^{2}) = 2  \Omega(X)X - \Omega(\textbf{1}) X^2 , \ for \ all \ X \in \mathcal{T}.
\end{align}
Comparing (2.6) and (2.7) implies that
\begin{align}
\Omega(X)X = X \Omega(X), \ for \ all \ X \in \mathcal{T}.
\end{align}
Defining the mapping $\Delta:\mathcal{T} \rightarrow \mathcal{T}$ by $\Delta(X)= \Omega(X) - \Omega(\textbf{1})X$ and using (2.8), it is observed that
\begin{align*}
\Delta(X^{2}) & = \Omega(X^2) - \Omega(\textbf{1})X^2 \\ & = 2  \Omega(X)X - \Omega(\textbf{1})X^2  - \Omega(\textbf{1})X^2 \\ & = \Omega(X)X - \Omega(\textbf{1})X^2  + X \Omega(X) - \Omega(\textbf{1})X^2 \\ & = (\Omega(X) - \Omega(\textbf{1})X)X + X(\Omega(X) - \Omega(\textbf{1})X) \\ & = \Delta(X)X + X \Delta(X)
\end{align*}
for all $X \in \mathcal{T}$, which means that $\Delta$ is a Jordan derivation on $\mathcal{T}$. Furthermore, since $\Omega(X)X = X \Omega(X)$ for all $X \in \mathcal{T}$, and $\Omega(\textbf{1}) \in Z(\mathcal{T})$, we infer that $\Delta(X)X = X \Delta(X)$ for all $X \in \mathcal{\mathcal{T}}$. Hence, $\Delta$ is a commuting Jordan derivation on $\mathcal{T}$. According to \cite[Theorem 2.4]{H-J}, $\Delta$ is zero and consequently, $\Omega = L_{\Omega(\textbf{1})}$. This equation with the assumption that $\Omega(\textbf{1}) \in Z(\mathcal{T})$ implies that $\Omega(XY) = \Omega(X) Y = X \Omega(Y)$ for all $X, Y \in \mathcal{T}$, which means that $\Omega$ is a two-sided centralizer on $\mathcal{T}$. Using this fact along with (3.3), we have \begin{align*}\Psi(X) & = \frac{ \gamma}{n} \Omega(X) + \frac{\gamma (n - 1)}{n} X \Omega(\textbf{1}) \\ & = \frac{ \gamma}{n} \Omega(X) + \frac{\gamma (n - 1)}{n}\Omega(X) \\  & = \gamma \Omega(X)\end{align*} for all $X \in \mathcal{T}$. Since $\Omega$ is a two-sided centralizer and we are assuming that $\gamma \in Z(\mathcal{T})$, $\Psi$ is a two-sided centralizer on $\mathcal{T}$, as desired.
\end{proof}

As consequences of Theorem \ref{1}, we obtain the following corollary.

\begin{corollary}\label{3}
Let $n > 1$ be an integer, let $\gamma$ be an invertible of $Z(\mathcal{T})$, and let $\Psi, \Omega:\mathcal{T} \rightarrow \mathcal{T}$ be additive mappings satisfying \begin{align*}
\Psi(X^n) = \gamma X\Omega(X^{n-1}) = \gamma \Omega(X^{n - 1}) X\end{align*}
for all $X \in \mathcal{T}$. If $\Omega(\textbf{1}) \in Z(\mathcal{T})$, then $\Psi$ and $\Omega$ are two-sided centralizers on $\mathcal{T}$ and also $\Psi = \gamma \Omega$.
\end{corollary}

\begin{proof} We proceed with the proof similarly to the proof of Theorem \ref{1}. According to the assumption, we have
\begin{align}
\Psi(X^n) = \gamma X\Omega(X^{n - 1}), \ \ \ \ \ \ \ \ \ X \in \mathcal{T}.
\end{align}
Let $C$ be an arbitrary element of $Z(\mathcal{T})$. Replacing $X$ by $X + C$ in (2.9), we have
\begin{align*}
\Psi\Big(\sum_{k = 0}^{n}\Big(^{n}_{k}\Big)X^{n - k}C^{k}\Big)  = \gamma (X + C)\Omega\Big(\sum_{k = 0}^{n-1}\Big(^{n-1}_{k}\Big)X^{n -1 - k}C^{k}\Big)
\end{align*}
Using equation (2.9) and collecting together terms of the above-mentioned expressions involving the same number of factors of $C$, it is obtained that
\begin{align}
\sum_{k = 1}^{n - 1}\Theta_k(X, C) = 0, \ \ \ \ \ \ \ \ \ \ \ \ \ X \in \mathcal{T},
\end{align}
where
\begin{align*}
\Theta_k(X, C) = \Big(^{n}_{k}\Big)\Psi(X^{n - k}C^{k}\Big) - \Big(^{n-1}_{k}\Big)\gamma X \Omega(C^{k}X^{n - k - 1}) - \Big(^{n-1}_{k-1}\Big) \gamma C \Omega(C^{k - 1}X^{n - k}).
\end{align*}
Replacing $C$ by $C, 2C, ..., (n - 1)C$ in equation (2.10), and using the method used in Theorem \ref{1}, we obtain that $\Theta_{n - 2}(X,C) = 0$. So
\begin{align*}
0 = & \Theta_{n - 2}(X,C) = \Big(_{n - 2}^{n}\Big)\Psi(X^{2}C^{n - 2}) - \Big(^{n-1}_{n-2}\Big) \gamma X \Omega(C^{n-2}X) -  \Big(^{n-1}_{n - 3}\Big) \gamma C \Omega(C^{n - 3}X^{2})
\end{align*}
Setting $C = \textbf{1}$ in the previous equality, we get that
\begin{align*}
\Big(_{n - 2}^{n}\Big)\Psi(X^{2}) =  \Big(^{n-1}_{n-2}\Big) \gamma X \Omega(X) + \Big(^{n-1}_{n - 3}\Big) \gamma \Omega(X^{2}), \ for \ all \ X \in \mathcal{T}.
\end{align*}
So, we can deduce that
\begin{align}
n \Psi(X^{2}) = 2 \gamma X \Omega(X) + (n-2) \gamma \Omega(X^{2}), \ for \ all \ X \in \mathcal{T}.
\end{align}
Considering the assumption that $\Psi(X^n) = \gamma \Omega(X^{n - 1}) X$ and using the above argument, we have
\begin{align}
n \Psi(X^{2}) = 2 \gamma \Omega(X)X + (n-2) \gamma \Omega(X^{2}), \ for \ all \ X \in \mathcal{T}.
\end{align}
In view of the equalities (2.11) and (2.12) and considering $\mu = n \Psi - (n -2)\gamma \Omega$, it is obtained that
\begin{align}
\mu(X^{2}) = 2 \gamma X \Omega(X) = 2 \gamma\Omega(X)X, \ for \ all \ X \in \mathcal{T}.
\end{align}

Using equation (2.13) and Theorem \ref{1}, we conclude that $\Omega$ is a two-sided centralizer and also $\mu = 2 \gamma \Omega$. So, we can write $n \Psi - (n -2)\gamma \Omega = \mu = 2 \gamma\Omega$, and this implies that $\Psi = \gamma \Omega$. Since $\Omega$ is a two-sided centralizer, so is $\Psi$, as required.
\end{proof}

In the next corollary, we offer a result concerning the automatic continuity of some linear mappings satisfying the functional identities mentioned in Theorem \ref{1} and Corollary \ref{3}.

\begin{corollary}\label{4}
Let $\mathcal{T}$ be a unital triangular normed algebra, let $n > 1$ be an integer, let $\gamma$ be an invertible of $Z(\mathcal{T})$ and let $\Psi, \Omega:\mathcal{T} \rightarrow \mathcal{T}$ be linear mappings such that $\Omega(\textbf{1}) \in Z(\mathcal{T})$. If $\Psi$ and $\Omega$ satisfy either of the following functional identities, then both are continuous.
\begin{enumerate}
	\item   [(i)] $\Psi(X^n) = \gamma X^{n - 1}\Omega(X) = \gamma \Omega(X) X^{n - 1}$, for all $X \in \mathcal{T}$;
	\item   [(ii)]  $\Psi(X^n) = \gamma X\Omega(X^{n-1}) = \gamma \Omega(X^{n - 1}) X$, for all $X \in \mathcal{T}$;
\end{enumerate}
\end{corollary}

\begin{proof} If the mappings $\Psi$ and $\Omega$ satisfy functional equations (i) or (ii), then it follows from Theorem \ref{1} and Corollary \ref{3} that $\Omega = L_{\Omega(\textbf{1})}$ and $\Psi = \gamma \Omega$. Evidently, both $\Omega$ and $\Psi$ are continuous.
\end{proof}

In the  theorem below, we present a functional equation characterizing two-sided generalized derivations.

\begin{theorem}\label{5}
Let $n > 1$ be an integer and let $\Psi, \Omega:\mathcal{T}\rightarrow \mathcal{T}$ be linear mappings satisfying \begin{align*}
2\Psi(X^n) =  X^{n - 1}\Omega(X) + \Omega(X) X^{n - 1}\end{align*}
for all $X \in \mathcal{T}$. If $\Omega(\textbf{1}) \in Z(\mathcal{T})$, then $\Psi$ and $\Omega$ are two-sided generalized derivations on $\mathcal{T}$.
\end{theorem}
\begin{proof}According to the aforementioned assumption, we have
\begin{align}
2\Psi(X^n) =  X^{n - 1}\Omega(X) + \Omega(X) X^{n - 1},\ for \ all \ X \in \mathcal{T}.
\end{align}
Let $C$ be an arbitrary element of $Z(\mathcal{T})$. Replacing $X$ by $X + C$ in (2.14), we have
\begin{align*}
2\Psi\Big(\sum_{k = 0}^{n}\Big(^{n}_{k}\Big)X^{n - k}C^{k}\Big)  = \sum_{k = 0}^{n - 1}\Big(^{n - 1}_{k}\Big)X^{n - k - 1}C^{k}\left(\Omega(X) + \Omega(C)\right) + \left(\Omega(X) + \Omega(C)\right)\sum_{k = 0}^{n - 1}\Big(^{n - 1}_{k}\Big)X^{n - k - 1}C^{k}
\end{align*}
Using equation (2.14) and collecting together terms of the above-mentioned expressions involving the same number of factors of $C$, it is obtained that
\begin{align}
\sum_{k = 1}^{n - 1}\Lambda_k(X, C) = 0, \ for \ all \ X \in \mathcal{T},
\end{align}
where
\begin{align*}
\Lambda_k(X, C) &= 2\Big(^{n}_{k}\Big)\Psi(X^{n - k}C^{k}\Big) - \Big(^{n-1}_{k-1}\Big)C^{k-1}X^{n - k}\Omega(C) - \Big(^{n-1}_{k-1}\Big)C^{k - 1}\Omega(C) X^{n - k} \\ & - \Big(^{n-1}_{k}\Big)C^{k}X^{n -1 - k}\Omega(X) - \Big(^{n-1}_{k}\Big)C^{k}\Omega(X)X^{n - 1- k}.
\end{align*}
Replacing $C$ by $C, 2C, ..., (n - 1)C$ in equation (2.15), and using the method used in Theorem \ref{1}, we get that $\Lambda_{n - 1}(X,C) = \Lambda_{n - 2}(X,C) = 0$. Since $\Lambda_{n - 1}(X,C) = 0$, we have
\begin{align*}
0 = & \Lambda_{n - 1}(X,C) = 2\Big(^{n}_{n-1}\Big)\Psi(X C^{n-1}\Big) - \Big(^{n-1}_{n-2}\Big)C^{n-2}X \Omega(C) - \Big(^{n-1}_{n-2}\Big)C^{n-2}\Omega(C)X - 2 C^{n-1}\Omega(X)
\end{align*}
for all $X \in \mathcal{T}$. Putting $C = \textbf{1}$ in the above equation, we get that
\begin{align}
\Psi(X) = \frac{ \Omega(X)}{n}  + \frac{n - 1}{n} X \Omega(\textbf{1}), \ for \ all \ X \in \mathcal{T}.
\end{align}
It follows from (2.16) that
\begin{align}
\Psi(X^2) = \frac{ \Omega(X^2)}{n}  + \frac{n - 1}{n} X^2 \Omega(\textbf{1}), \ for \ all \ X \in \mathcal{T}.
\end{align}

On the other hand, since $\Lambda_{n - 2}(X,C) = 0$, we have the following expressions:
\begin{align*}
0 = & \Lambda_{n - 2}(X,C) = 2\Big(^{n}_{n-2}\Big)\Psi(X^{2}C^{n-2}\Big) - \Big(^{n-1}_{n-3}\Big)C^{n-3}X^{2}\Omega(C) - \Big(^{n-1}_{n-3}\Big)C^{n-3}\Omega(C) X^{2} \\ & - \Big(^{n-1}_{n-2}\Big)C^{n-2}X \Omega(X) - \Big(^{n-1}_{n-2}\Big)C^{n-2}\Omega(X)X
\end{align*}
for all $X \in \mathcal{T}$. Putting $C = \textbf{1}$ in the above equation, we get that
\begin{align}
\Psi(X^{2}) = \frac{X \Omega(X)}{n} + \frac{\Omega(X)X}{n} + \frac{n - 2}{n} X^{2} \Omega(\textbf{1}), \ for \ all \ X \in \mathcal{T}.
\end{align}
Comparing (2.17) and (2.18), we get that
\begin{align*}
\Omega(X^{2}) = X \Omega(X) + \Omega(X)X - \Omega(\textbf{1}) X^2, \ for \ all \ X \in \mathcal{T}.
\end{align*}
Defining the mapping $\Delta: \mathcal{T} \rightarrow \mathcal{T}$ by $\Delta(X) = \Omega(X) - \Omega(\textbf{1})X$, it is observed that $\Omega(X^2) = \Omega(X)X + X\Delta(X) = \Delta(X)X + X\Omega(X)$ for all $X \in \mathcal{T}$. Also, note that
\begin{align*}
\Delta(X^2) & = \Omega(X^2) - \Omega(\textbf{1})X^2 \\ & = X \Omega(X) + \Omega(X)X - \Omega(\textbf{1}) X^2 - \Omega(\textbf{1})X^2 \\ & = (\Omega(X) - \Omega(\textbf{1}) X)X + X (\Omega(X) - \Omega(\textbf{1}) X) \\ & = \Delta(X)X + X \Delta(X)
\end{align*}
for all $X \in \mathcal{T}$, which means that $\Delta$ is a Jordan derivation on $\mathcal{T}$. According to the main theorem of \cite{F}, $\Delta$ is a derivation on $\mathcal{T}$. So, we have the following expressions:
\begin{align*}
\Omega(XY) &= \Delta(XY) + \Omega(\textbf{1})XY \\ & = \Delta(X)Y + X\Delta(Y) + \Omega(\textbf{1})XY \\ & = (\Delta(X) + \Omega(\textbf{1})X)Y + X \Delta(Y) \\ & = \Omega(X)Y + X \Delta(Y)
\end{align*}
for all $X, Y \in \mathcal{T}$. Easily, one can see that $$\Omega(XY)= \Delta(X)Y + X \Omega(Y)$$ for all $X, Y \in \mathcal{T}$, and this means that $\Omega$ is a two-sided generalized derivation on $\mathcal{T}$. Our next task is to show that $\Psi$ is also a two-sided generalized derivation on $\mathcal{T}$. In view of (2.16), we have
\begin{align*}
\Psi(X) = \frac{ \Omega(X)}{n}  + \frac{n - 1}{n} X \Omega(\textbf{1}) = \frac{ \Delta(X) + \Omega(\textbf{1})X}{n}  + \frac{n - 1}{n} X \Omega(\textbf{1}) = \frac{\Delta(X)}{n} + \Omega(\textbf{1})X,
\end{align*}
for all $X \in \mathcal{T}$. A straightforward verification shows that $\Psi$ is a two-sided generalized derivation.
\end{proof}

Here, a direct consequence of the last theorem can be stated. This result is considered as an extension of \cite[Theorem 1]{V7}.

\begin{corollary}
Let $n > 1$ be an integer and let $\Psi:\mathcal{T}\rightarrow \mathcal{T}$ be linear mappings satisfying \begin{align*}
2\Psi(X^n) =  X^{n - 1}\Psi(X) + \Psi(X) X^{n - 1}\end{align*}
for all $X \in \mathcal{T}$. If $\Psi(\textbf{1}) \in Z(\mathcal{T})$, then $\Psi$ is a two-sided centralizer.
\end{corollary}

\begin{proof}
We can assume in Theorem \ref{5} that $\Psi = \Omega$. Based on what was presented in the proof of Theorem \ref{5}, it is observed that $\Delta + L_{\Omega(\textbf{1})} = \Omega = \Psi = \frac{\Delta}{n} + L_{\Omega(\textbf{1})}$. Since $n > 1$, we conclude that $\Delta = 0$, and so $\Psi$ is a two-sided centralizer, as desired.
\end{proof}


\bibliographystyle{amsplain}

\vskip 0.5 true cm 	

{\tiny (Amin Hosseini) Kashmar Higher Education Institute, Kashmar, Iran}
	
	{\tiny\textit{E-mail address:} a.hosseini@kashmar.ac.ir}

\vskip 0.3 true cm 	
\end{document}